\documentclass[10pt,twoside]{article}
\usepackage[psamsfonts]{amssymb}
\usepackage{amsfonts}
\usepackage[mathcal]{eucal}
\usepackage{Latex-document}

\newcommand{\g}{\mathfrak{g}}
\renewcommand{\a}{\mathfrak{a}}

\newcommand{\q}{\mathfrak{q}}

\newcommand{\n}{\mathfrak{n}}

\newcommand{\s}{\mathfrak{s}}

\newcommand{\N}{\mathbb{N}}
\newcommand{\F}{\mathbb{F}}

\newcommand{\R}{\mathbb{R}}
\newcommand{\C}{\mathbb{C}}
\newcommand{\D}{\mathbb{D}}

\newcommand{\qed}{{\null\hfill\
\raise3pt\hbox{\framebox[0.1in]{}}}\break\null}

\newcommand{\ste}{\hfill\break}

\title{\bf Harmonic Analysis on Real Reductive \vskip -2mm Symmetric
Spaces\vskip 6mm}
\author{Patrick Delorme\vspace*{-0.5cm}\thanks{Institut
de Math\'ematiques de Luminy, U.P.R. 9016 du C.N.R.S., Facult\'e
des Sciences de Luminy, 163 Avenue de Luminy, Case 930, 13288
Marseille Cedex 09, France. E-mail: delorme@iml.univ-mrs.fr}}
\date{\vspace{-8mm}}

\begin{document}

\maketitle

\thispagestyle{first} \setcounter{page}{545}

\begin{abstract}

\vskip 3mm

Let $G$ be a reductive group in the Harish-Chandra class e.g. a
connected semisimple Lie group with finite center, or the group of
real points of a connected reductive algebraic group defined over
$\R$. Let $\sigma$ be  an involution of the Lie group $G$, $H$ an
open subgroup of the subgroup of fixed  points of $\sigma$. One
decomposes the elements of $L^2(G/H)$ with the help of joint
eigenfunctions under the algebra of left invariant differential
operators under $G$ on $G/H$.

\vskip 4.5mm


\end{abstract}

\vskip 12mm

\section{Introduction} \setzero

\vskip -5mm \hspace{5mm}

Let $G$ be a real reductive group in the Harish-Chandra class
[H-C1],  e.g. a connected semisimple Lie group with finite center,
or  the group of real points of a connected reductive algebraic
group defined over $\R$. Let $\sigma$ be  an involution of the Lie
group $G$, $H$ an open subgroup of the subgroup of fixed points of
$\sigma$.

Important problems of harmonic analysis on the so-called
reductive symmetric space $G/H$ are  :

{(a) to make the simultaneous spectral decomposition of the
elements of the algebra $\D (G/H)$ of left invariant differential
operators under $G$ on $G/H$. In other words, one wants to write
the elements of $L^2(G/H)$ with the help of joint eigenfunctions
under $\D (G/H)$.}

{(b) to decompose the left regular representation of $G$ in
$L^2(G/H)$ into an Hilbert integral of irreducible unitary
representations of $G$ : this is essentially the Plancherel
formula.}

{(c) to decompose the Dirac measure at $eH$, where $e$ is the
neutral element of $G$, into an integral of $H$-fixed
distribution vectors : this is essentially the Fourier inversion
formula.}

These problems were solved for the ``group case" (i.e. the group
viewed as a symmetric space : $G=G_{1}\times G_{1}$,
$\sigma(x,y)=(y,x)$, $H$ is the diagonal of $G_{1}\times G_{1}$ )
by Harish-Chandra in 1970s (see [H-C1,2, 3]), the Riemannian case
($H$ maximal compact) had been treated before (see [He]). Later,
there were deep results by T. Oshima [O1]. When $G$ is a complex
group and $H$ is a real form, the Problems (a), (b), (c) were
solved by P. Harinck, together with an inversion formula  for
orbital integrals ([Ha], see also [D3] for the link of her work
with the work of A. Bouaziz on real reductive groups).

Then,  E. van den Ban and  H. Schlichtkrull, on the one hand, and
I, on the other hand, obtained different  solutions to problems
(a), (b), (c). Moreover, they obtained a Paley-Wiener theorem
(see [BS3] for a presentation of  their work). I present here my
point of view, with an    emphasize on problem (a), because it
simplifies  the formulations of the results (nevertherless, the
important aspect of representation theory is hidden). It includes
several joint works, mainly with J. Carmona , and also with E.
van den Ban and J.L. Brylinski. Severals works  of T. Oshima,
linked to the the Flensted-Jensen duality, alone and with  T.
Matsuki are  very important in my proof, as well as earlier
results of E. van den Ban and H. Schlichtkrull.

I have to acknowledge the deep  influence of Harish-Chandra 's
work. The crucial role played by the work [Be] of J. Bernstein on
the support of the Plancherel measure, and some part of Arthur's
article on the local trace formula [A] will be appearant in the
main body of the article.

\section{Temperedness of the spectrum} \setzero

\vskip -5mm \hspace{5mm}

Let $\theta$ be a Cartan involution of $G$ commuting with
$\sigma$, let $K$ be the fixed point set of $\theta$. Let $\g$ be
the Lie algebra of $G$, etc. Let $\s$ (resp. $\q$) be the space
of elements in $\g$ which are antiinvariant under the
differential of $\theta$ (resp. $\sigma$). Let $\a_{\emptyset}$
be a maximal abelian subspace of $\s\cap \q$. If $P$ is a
$\sigma\theta$-stable parabolic subgroup of $G$, containing
$A_{\emptyset}:=exp\>\a_{\emptyset}$, we denote by
$P=M_{P}A_{P}N_{P}$ its Langlands $\sigma$-decomposition. More
precisely  $A_{P}$ is the subgroup of the elements $a$ of the
split component of the Levi factor $L_{P}=P\cap\theta(P)$ such
that $\sigma(a)=a^{-1}$. Here $M_{P}$ is larger than that for the
usual Langlands decomposition.

In order to simplify the exposition we will make the following
simplifying assumption:

\noindent {\bf Hypothesis: } {\it For any $P$ as above, $HP$ is
the unique  open $(H,P)$-double coset.}

When $\sigma=\theta$ (the case of a riemannian symmetric space )
or the ``group case", this hypothesis is satisfied.

To get the Plancherel formula, it is useful to use $K$-finite
functions. They are often replaced by $\tau$-spherical functions.
Here $(\tau,V_{\tau})$ is a finite dimensional unitary
representation of $K$ and a $\tau$-spherical function on $G/H$ is
a function $ f:G/H\rightarrow V_{\tau}$ such that
$f(kx)=\tau(k)f(x)$, $k\in K$, $x\in G/H$.

Some spaces of $\tau$-spherical functions on $G/H$ play a crucial
role in the theory, namely:

(a) ${\cal C}(G/H,\>\tau)$ : the Schwartz space of
$\tau$-spherical functions on $G/H$ which are rapidly decreasing
as well as their derivatives by elements of the enveloping
algebra $U(\g)$ of $\g$ (see [B2]).

(b) ${\cal A}(G/H,\>\tau)$: the space of smooth  $\tau$-spherical
functions on $G/H$ which are $\D (G/H)$ finite. Here ${\cal A}$
is used to evoke automorphic forms.

(c) ${\cal A}_{temp}(G/H,\>\tau)$: the space of elements of ${\cal
A}(G/H,\>\tau)$ which have tempered growth as well as their
derivatives by elements of  $U(\g)$ ([D2]). Integration of
functions on $G/H$ defines a pairing between ${\cal
A}_{temp}(G/H,\>\tau)$ and ${\cal C}(G/H,\>\tau)$.

(d) ${\cal A}_{2}(G/H,\>\tau)$: the space of square integrable
elements of ${\cal A}(G/H,\>\tau)$. This is a subspace of the
three preceeding spaces.

One has:

\noindent {\bf Theorem 1} ( [D2] ){\bf: } {\it The space ${\cal
A}_{2}(G/H,\>\tau)$ is finite dimensional.}

This is deduced from the theory of discrete series for $G/H$
initiated by M. Flensted-Jensen [F-J] and  achieved by T.
Oshima   and T. Matsuki, using the Flensted-Jensen duality [OM].
One has also to use the behaviour of  the discrete series under
certain translation functors, studied by D. Vogan [V] and a
result of H. Schlichtkrull [S] on the minimal $K$-types of
certain discrete series.

The next result follows from the work of J. Bernstein [Be] on the
support of the
Plancherel measure.

\noindent {\bf Theorem 2} ([CD1], Appendice C){\bf : } {\it Every
function in ${\cal C}(G/H,\>\tau)$ can be canonically
desintegrated as an integral of elements of  ${\cal
A}_{temp}(G/H,\>\tau)$.}

This information appeared to be crucial at the end of our proof.

\section {The continuous spectrum: Eisenstein integrals} \setzero

\vskip -5mm \hspace{5mm}

Let $P=MAN$ the Langlands $\sigma$-decomposition of a
$\sigma\theta$-stable parabolic subgroup $P$ of $G$. Let
$\rho_{P}$ be the half sum of the roots of $\a$ in $\n$ and
$\lambda\in \a^{*}_{\C}$ be such that the real part of $\lambda
 -\rho_{P}$ is strictly dominant with respect to the roots
of $\a$ in $\n$. Let $\tau_{M}$ be the restriction of $\tau$ to
$M\cap K$. Then, if $x\in G/H$ and $\psi\in {\cal A}_{2}(M/M\cap
H,\tau_{M})$, the following integral is convergent:
$$E(P,\psi,\lambda)(x):=\int _{K}\tau
(k^{-1})\Psi_{\lambda}(kx)dk,$$ where $\Psi_{\lambda}(x)=0$ if
$x\notin PH$, and $\Psi_{\lambda}(x)=a^{-\lambda
+\rho_{P}}\psi(m)$ if $x=namH$ with $n\in N$, $a\in A$, $m\in M$.
Moreover $E(P,\psi, \lambda)$ is an element of ${\cal
A}(G/H,\>\tau)$. Eisenstein integrals are the $\tau$-spherical
versions of $K$-finite functions of the form: $gH\mapsto<\pi
'(g)\xi,v>$, where $\pi'$ is the contragredient of a generalized
principal series $\pi$, $\xi$ is a certain $H$-fixed distribution
vector of $\pi$, $v$ is a $K$-finite vector of $\pi$.

\noindent{\bf Theorem 3} ([BrD]){\bf : } {\it The function
$\lambda\mapsto E(P,\psi,\lambda)$ admits a meromorphic
continuation in $\lambda\in \a^{*}_{\C}$. This meromorphic
continuation, denoted in the same way, multiplied by a suitable
product, $p_{\psi}$, of functions of type $\lambda\mapsto
(\alpha, \lambda)+c$, where $\alpha$ is a root of $\a$ and $c\in
\C$, is holomorphic around $i\a^{*}$.}

This meromorphic continuation is an interesting feature of the
theory. For the the ``group case", it comes down to the
meromorphic continuation of Knapp-Stein intertwining integrals.
My  proof with Brylinski uses $D$-modules arguments.

The case where $P$ is minimal had been treated separately by E.
van den Ban [B1] and G. Olafsson [Ol]. One has also to mention the
earlier work of T. Oshima and J. Sekiguchi [OSe] on the spaces of
type $G/K_{\varepsilon}$.

The proof which gives the best results uses a method of tensoring
by finite dimensional representations. It is a joint work with
J.  Carmona. It was initiated by D. Vogan and N. Wallach (see
[W], chapter 10) for the meromorphic continuation of the
Knapp-Stein intertwining integrals. For symmetric spaces and the
most continuous spectrum, the  proof is due to E. van den Ban
[B2]. This proof uses Bruhat's thesis and tensoring by finite
dimensional modules. This implies rough estimates for Eisenstein
integrals, which generalize those obtained  by E. van den Ban
when $P$ is minimal [B2].

\noindent {\bf Theorem 4} ([D1]){\bf : } {\it If $\lambda\in
i\a^{*}$ is such that $E(P,\psi,\lambda)$ is well defined, then
it is tempered, i.e. is an element of ${\cal
A}_{temp}(G/H,\tau)$}.

This is a natural result but the proof is quite long. It uses the
behaviour under translation functors of $H$-fixed distribution
vectors of discrete series and of generalized principal series,
and also of the Poisson transform. Moreover the duality of M.
Flensted-Jensen, [F-J], and a criteria of temperedness due to
Oshima [O2] play a crucial role (apparently, J. Carmona has a way
to avoid boundary values).

With the help of this theorem and by using techniques due to E.
van den Ban [B2], the rough estimates for Eisenstein integrals
can be improved to get uniform sharp estimates for
$p_{\psi}(\lambda)E(P,\psi,\lambda)$, $\lambda\in i\a^{*}$ (cf.
[D1]).

\section{\boldmath$C$-functions} \setzero

\vskip -5mm \hspace{5mm}

Let $P$ be as above and let $L$ be equal to $MA$. The theory of
the constant term, due to J. Carmona [C1] (Harish-Chandra for the
group case, [H-C1]), gives a linear map from ${\cal
A}_{temp}(G/H,\>\tau)$ into ${\cal A}_{temp}(L/L\cap
H,\>\tau_{L})$, $\varphi \mapsto \varphi_{P}$, characterized by :
$$lim_{t\rightarrow +\infty}\> \delta_{P}^{1/2}((exp\>
tX)l)\varphi((exp\> tX)l) -\varphi_{P}((exp\> tX)l)=0,$$ where
$l\in L/L\cap H$, $X\in \a_{P}$ is $P$-dominant and $\delta_{P}$
is the modular function of $P$.

Let $Q$ be a $\sigma\theta $-stable parabolic subgroup of $G$
with the same $\theta$-stable Levi subgroup $L$ other than $P$.
Let $W(\a)$ be the group of automorphisms of $\a$ induced by an
element of $Ad(G)$. One defines meromorphic functions on
$\a^{*}_{\C}$, $\lambda \mapsto C_{Q\mid P}(s, \lambda)$ with
values in $End({\cal A}_{2}(M/M\cap H,\tau_{M}))$ such that :
$$E(P,\psi, \lambda)_{Q}(ma)=\sum _{s\in W(\a)}(C_{Q\mid P}(s,
\lambda)\psi)(m)a^{-s\lambda},\> m\in M,\>a\in A, \> \lambda \in
i\a^{*},$$ or rather for $\lambda$ in an open dense subset of $
i\a^{*}$.

The $C$-functions allow to normalize Eisenstein integrals as
folllows:
$$E^{0}(P,\psi, \lambda):=E(P,C_{P\mid P}(1,\lambda)^{-1}\psi
,\lambda).$$

\section{Truncation, Maass-Selberg relations and the
regularity of normalized Eisenstein integrals}  \setzero

\vskip -5mm \hspace{5mm}

Let $P$ be as above and let $P'=M'A'N'$ be the Langlands
$\sigma$-decomposition of another $\sigma\theta$-stable parabolic
subgroup of $G$. Let $\psi$ (resp. $\psi '$) be an element of
${\cal A}_{2}(M/M\cap H,\tau_{M})$ (resp. ${\cal A}_{2} (M'/M'\cap
H,\tau_{M'})$). One chooses $p_{\psi}$  as in Theorem 3, such that
the product of $p_{\psi}$ with the $C$-functions are holomorphic
in a neighbourhood of $i\a^{*}$ which is a product of $i\a^{*}$
with a neighbourhood of 0 in $\a^{*}$ . We do the same for $\psi
'$. One defines $F(\lambda):=p_{\psi}(\lambda)E(P,\psi, \lambda)$.
One defines similarly $F'$.

One assumes, for the rest of the article, that $G$ is semisimple.
This is just to simplify the exposition. One chooses $T\in
\a_{\emptyset}$, regular with respect to the roots of
$\a_{\emptyset}$ in $\g$. Let $C^{1}_{T}$ be the convex hull of
$W(\a_{\emptyset})T$ and let $C_{T}$ be equal to the subset
$K(exp\> C^{1}_{T})H$ of $G/H$.

\noindent{\bf Theorem 5} ([D2]){\bf : }

(i) {\it One gets an explicit expression $\omega^{T}(\lambda,
\lambda ')$, involving the $C$-functions (see an example below)
and vanishing when $A$ and $A'$ are not conjugate under $K$,
which is asymptotic to
$$\Omega^{T}(\lambda, \lambda '):=\int
_{C_{T}}(F(\lambda)(x),F'(\lambda ')(x))dx,$$
when $T$ goes to
infinity and $\lambda\in i \a^{*}$, $\lambda '\in i \a '^{*}$.
More precisely for $\delta
>0$ there exists $C>0$, $k\in \N$ and $\varepsilon >0$,
such that for all $T$ satisfying $\mid \alpha (T)\mid \geq \delta
\mid\mid T\mid\mid$ for every root $\alpha$ of $\a_{\emptyset}$
in  $\g$, one has:
$$\mid\Omega^{T}(\lambda,
\lambda ')-\omega^{T}(\lambda, \lambda ')\mid\leq C(1+ \mid
\mid\lambda\mid\mid)^{k} (1+ \mid \mid\lambda
'\mid\mid)^{k}e^{-\varepsilon \mid\mid T\mid\mid}.$$}

(ii) {\it Moreover $\omega^{T}$ is analytic in $(\lambda, \lambda
')\in i\a^{*}\times i\a '^{*}$.}

\noindent This generalizes a result of J. Arthur for the group
case [A], Theorem 8.1. My proof is quite similar, but I was able
to avoid his use of the Plancherel formula.

(ii) is an easy consequence of (i). In fact, the explicit form of
$\omega^T$ implies that it is meromorphic around $i\a^{*}\times
i\a '^{*}$. Moreover $\Omega^T$ is holomorphic, hence  locally
bounded, around $i\a^{*}\times i\a '^{*}$. From the inequality in
(i), one deduces that $\omega ^T$ is locally bounded, hence
holomorphic, around $i\a^{*}\times i\a '^{*}$.

We will now show, by an example, how the explicit form of
$\omega^{T}$ and its analyticity in $(\lambda,\lambda ')\in
i\a^{*}\times i\a '^{*}$ imply the Maass-Selberg relations.

Let $\sigma$ be equal to $\theta$, $H=K$, and $\tau$ be the
trivial representation. Let $P,P'$ be minimal parabolic subgroups
of $G$. Then $dim {\cal A}_{2}(M/M\cap H, \tau)=1$ and the
$C$-functions are scalar valued. One assumes $\g$ to be
semisimple and that the dimension of $A$ is one. Then $W(\a )$
has 2 elements, $\pm 1$, and one has the following explicit
expression of $\omega^{T}$:
$$\omega^{T}(\lambda, \lambda ')
=p_{\psi}(\lambda)p_{\psi '}(\lambda ')\sum_{s=\pm 1, s'=\pm
1}e^{s\lambda T- s'\lambda 'T}C_{P\mid P}(s,\lambda)\overline
{C_{P\mid P}(s',\lambda ')}(s\lambda-s' \lambda ')^{-1}.$$
Thus
$\omega^{T}(\lambda, \lambda ')$ is the sum of a product of
$(\lambda -\lambda ')^{-1}$ by an analytic function with a
product of $(\lambda + \lambda ')^{-1}$ by an analytic function.
The analyticity at $(\lambda, \lambda)$ implies easily that the
factor in front of  $(\lambda
 -\lambda ')^{-1}$  vanishes for $\lambda= \lambda '$.
Hence we get  $\mid C_{P\mid P}(1,\lambda)\mid ^2=\mid C_{P\mid
P}(-1,\lambda)\mid ^2$, $\lambda\in i\a^{*}$. This is one of the
Maass-Selberg relations(cf. [D2], Theorem 2, and the work with J.
Carmona [CD2], Theorem 2 for the general case, see [B1], [B2] for
the case where $P$ is minimal). These relations imply that the
$C$-functions attached to normalized Eisenstein integrals are
unitary, when defined, for $\lambda$ purely imaginary. Hence they
are locally bounded . This implies that they are holomorphic
around the imaginary axis. This implies in particular some
holomorphy property of the constant term of normalized Eisentein
integrals. From this, with the help of [BCD], one deduces:

\noindent {\bf Theorem 6} (Regularity theorem for normalized
Eisenstein integrals, [CD2], [BS1]  for $P$ minimal){\bf : } {\it
The normalized Eisenstein integrals are holomorphic in a
neighbourhood of the imaginary axis}.

\section{Fourier transform and wave packets} \setzero

\vskip -5mm \hspace{5mm}

\noindent {\bf Theorem 7} ([CD2],  [BS1] for {\boldmath $P$}
minimal){\bf : } {\it For $f\in {\cal C}(G/H,\tau)$, one has
${\cal F}^0_{P}f\in {\cal S}(i\a^{*})\otimes {\cal A}_2 (M/M\cap
H,\tau _M)$, where ${\cal F}^0_{P}f$ is characterized by:
$$(({\cal F}^0_{P}f)(\lambda),
\psi)=\int _{G/H} (f(x),E^0 (P,\psi,\lambda)(x))dx, \> \> \>
\>\lambda\in i\a^*,\> \>\> \>\psi\in  {\cal A}_2 (M/M\cap H, \tau
_M),$$ here ${\cal S}(i\a^{*})$ is the usual Schwartz space.}

This theorem follows from the sharp estimates of Eisenstein
integrals.

\noindent {\bf Theorem 8} ([BCD]){\bf : } {\it If $\Psi$ is an
element of ${\cal S}(i\a^{*})\otimes{\cal A}_2 (M/M\cap H, \tau
_M)$, one has ${\cal I} ^{0}_{P}\in {\cal C} (G/H,\tau)$, where :
$$ {\cal I}
^{0}_{P}(x):=\int _{i\a ^{*}}E^0 (P,\Psi(\lambda),
\lambda)d\lambda,\>\> \>\>x\in G/H.$$}

This theorem follows from the regularity theorem and from the
joint work with E. van den Ban and J. Carmona, [BCD].

Now we want to compute ${\cal F}^{0}_{P'}{\cal I}^{0}_{P}$. For
this purpose one has to study the integral:
$$I:=\int _{G/H}(\int _{i\a^{*}}
\alpha(\lambda)E^{0}(P,\psi,\lambda)(x) d\lambda,E^{0}(P',\psi
',\lambda ')(x))dx.$$

One truncates the integral on $G/H$ to $C_{T}$ and let $T$ goes
to infinity (far from the walls). Let us denote the truncated
inner product of the normalized Eisenstein integrals by
$^0\Omega^T(\lambda,\lambda ')$. Using Fubini's theorem one has:
$$I=lim_{T\rightarrow\infty}\int_{i\a^{*}}
\alpha(\lambda)^0\Omega^T(\lambda,\lambda ')d\lambda.$$ As for unnormalized Eisentein integrals, one has an
asymptotic evaluation of \linebreak $^0\Omega^T(\lambda,\lambda ')$ by an explicit expression
$^0\omega^T(\lambda,\lambda ')$. One can replace $^0\Omega^T$ by $^0\omega^T$ in the previous formula. By using an
expression of $^0\omega^T(\lambda,\lambda ')$, viewed as a distribution in $\lambda$, for fixed $\lambda '$ ,
involving Fourier transforms of cones ( [D2], Theorem 3) one gets:

\noindent {\bf Theorem 9} ([CD2]){\bf : } {\it Let $\F$ be a set
of representative of $\sigma$-association classes of
$\sigma\theta$-stable parabolic subgroups. Here
$\sigma$-association means that the $\a$ are conjugated under $K$.
Define:
$${\cal P}_{\tau}=\sum_{P\in \F}((W(\a_{P}))^{-1}
{\cal I}^0_P{\cal F}^0_P.$$ Then ${\cal P}_{\tau}$ is an
orthogonal projection operator in ${\cal C}(G/H,\tau)$ endowed
with the $L^2$ scalar product.}

\section{The Plancherel formula} \setzero

\vskip -5mm \hspace{5mm}

Essentially, the solution to problem (a) is contained in the
following:

\noindent {\bf Theorem 10} ([D4]){\bf : } {\it The projection
${\cal P}_{\tau}$ is the identity operator on ${\cal
C}(G/H,\tau)$.}

Actually, this gives an expression of every element in ${\cal
C}(G/H,\tau)$ as a wave packet of normalized Eisenstein integrals.
The proof goes as follows. If ${\cal P}_{\tau}$ was not the
identity, using Theorem 1 on the temperedness of the spectrum,
one could find a non zero element of ${\cal A}_{temp}(G/H,\tau)$
which is orthogonal to the image of ${\cal P}_{\tau}$. Then,
generalizing Theorem 5  to the truncated inner product of an
Eisenstein integral with a general element of ${\cal
A}_{temp}(G/H,\tau)$, this orthogonality can be explicitely
described (cf. the evaluation of $I$ before Theorem 9). As a
result, this function has to be zero, a contradiction which
proves the theorem.

The theorem translates to $K$-finite functions, involving
representations and $H$-fixed distribution vectors.

The theorem can also be expressed with the unnormalized
Eisenstein integrals. Then there are certain Plancherel factors
involved. They are linked, as in the group case , to the
intertwining integrals. Following the approach of A.Knapp and
G.Zuckerman, [KZ], their computation is reduced to find an
embedding of discrete series into principal series attached to
minimal parabolic subgroups. For connected groups, this has been
done by J. Carmona [C2].

\section{Applications and open problems} \setzero

\vskip -5mm \hspace{5mm}

\noindent {\bf Schwartz space for the hypergeometric Fourier
transform}

The image of a natural  Schwartz space by the hypergeometric
Fourier transform is characterized [D5]. The work uses the
Plancherel formula of E. Opdam [Op], and the techniques mentioned
above : theory of the constant term, $C$-functions, truncation
...

\noindent {\bf Generalized Schur orthogonality relations }

Using the Plancherel formula for reductive symmetric spaces
groups, K. Ankabout, [An], has proved generalized Schur
orthogonality relations for generalized coefficients related to
real reductive symmetric spaces. In particular, at least if we
assume multiplicity one in the Plancherel formula, it implies the
following:

There exists an explicit positive function $d$ on $G/H$, such
that,  for almost all representations $\pi$ occuring in the
Plancherel formula, for $\xi$ an $H$-fixed distribution vector of
$\pi$ occuring in the decomposition of the Dirac measure, there
exists an explicit non zero constant $C_{\pi}$ such that, for all
$v, v'$, $K$-finite vectors in the space of $\pi$:
$$\lim_{\varepsilon \rightarrow 0^{+}}\varepsilon
^{n_{\pi}}\int_{G/H}e^{-\varepsilon d(x)}<\pi'(g)\xi,v> {\overline
{<\pi'(g)\xi,v'>}} dx= C_{\pi}(v, v').$$

Here $n_{\pi}$ is the dimension of the support of the Plancherel
decomposition, around $\pi$. This refines and generalizes a work
of Mirodikawa. It suggests to look for such
type of relations in  other situations.\\

\noindent{\bf {\boldmath $\D(G/H)$}-finite {\boldmath
$\tau$}-spherical functions on reductive symmetric spaces }

S. Souaifi [So] showed how these functions appear as linear
combinations of derivatives along the complex parameter
$\lambda$, of Eisenstein integrals. For $K$-finite functions,
filtrations are introduced, whose subquotients are described in
terms of induced representations. The starting point is an
adaptation of ideas used by J. Franke to study spaces of
automorphic forms. The use of the spectral decomposition by
Langlands is replaced here by the use of the Plancherel formula.
For reductive $p$-adic groups, and for the group case, I got
similar results.

\noindent{\bf Invariant harmonic analysis on real  reductive
symmetric spaces }

The goal is to study the $H$-invariant eigendistributions under
$\D(G/H)$ on $G/H$ and to express invariant measures on certain
$H$-orbits  in terms of these distributions ( cf [D3] for the
work of A.Bouaziz and P.Harinck for the group case and
$G(\C)/G(\R)$, see also [OSe]).

\noindent {\bf Harmonic analysis on {\boldmath $p$}-adic reductive
symmetric spaces}

For the group case, the Problems (b) and (c) of the Introduction
have been solved by Harish-Chandra, up to the explicit
description of the discrete series. In general, the problems are
open (see [HH] for interesting structural results).

\label{lastpage}

\end{document}